# The universality of the Kalman filter: a conditional characteristic function perspective

Sandhya Rathore[1], Shambhu N. Sharma[2] and Shaival H. Nagarsheth[3]
Electrical Engineering Department,
S.V. National Institute of Technology, Surat, 395007, Gujarat, India.
E-mail: [1]sandhya.rathore.svnit@gmail.com, [2]snsvolterra@gmail.com
[3]shn411@gmail.com

**Abstract**

The universality of the celebrated Kalman filtering can be found in control theory. The Kalman filter has found its striking applications in sophisticated autonomous systems and smart products, which are attributed to its realization in a single complex chip. In this paper, we revisit the Kalman filter from the perspective of 'conditional characteristic function evolution and Itô calculus' and develop three Kalman filtering Theorems and their formal proof. Most notably, this paper reveals the following: (i) Kalman filtering equations are a consequence of the 'evolution of conditional characteristic function' for the linear stochastic differential system coupled with the linear discrete measurement system. (ii) The Kalman filtering is a consequence of the 'stochastic evolution of conditional characteristic function' for the linear stochastic differential system coupled with the linear continuous measurement system. (iii) The structure of the Kalman filter remains invariant under two popular stochastic interpretations, the Itô vs Stratonovich.

**Introduction**

Norbert Wiener in 1921 proved that the Brownian motion has continuous sample path and non-differentiable with probability one [1]. Norbert Wiener developed the filtering equation using the convolution integral equation for linear observation equations. The Wiener filter is a signal estimation method. Suppose the Wiener filter is the LTI system, the Wiener filtering can be regarded as an infinite-dimensional LTI system, the transfer function would be irrational. In 1960, Kalman sowed the seed of the filtering of dynamic systems in which he considered linear stochastic differential equation coupled with linear observation equation [2]. The result was published in the *Journal of Basic Engineering*. This attracted attention in greater circles via its Apollo mission application. The Kalman filtering has found appealing applications in smartphones as well. Thus, the Kalman filtering is ubiquitous and is in everybody's life. Suppose the Kalman filter is an LTI system, where the observation is the input and the filtered estimate is the output, that can be regarded as the finite-dimensional and the transfer function of the Kalman filter would be rational. The Kalman filtering and its discoverer's legacy can be found in appealing compilations and papers, e.g. [3]. Kalman meets Shannon [4] in noisy Gaussian channels, the ubiquitous Riccati equation is a part of the coupled Kalman filtering equations for the continuous state-continuous measurement systems. That reveals the conceptual depth of the Kalman filter. After combining the celebrated Wiener-Höpf equation, Wiener filter, linear stochastic differential equation, linear measurement equation coupled with the orthogonal projection Lemma, we arrive at the Kalman gain and coupled filtering equations. Furthermore, the notion of non-linear Kalman filtering was introduced. The non-linear filtering has greater conceptual depth as well as involves greater mathematical rigour. Non-linear filtering theory hinges on the filtering density evolution equation and the stochastic evolution of the conditional expectation of the scalar function. A remarkable paper on non-linear filtering can be traced back to a pioneering work of 1967, H J Kushner. That is popularly known the Kushner-Stratonovich equation, see [5]. An alternative interpretation of the Kushner-Stratonovich is FKK filtering [6]. In this paper, we restrict our discussions to the celebrated Kalman filtering.

This paper is inspired from the following two papers: (i) A simple concept with new perspectives leads to surprising results ([7]) (ii) Moore and Anderson published an interesting paper [8], the Kalman filtering: whence, what and whither. In their paper, they explained 'the Kalman filtering by extending it in non-linear filtering perspective with some discussions ([8], p. 52). In this paper, we answer the question, *whither of the Kalman filtering?* in greater detail. This paper sketches the answer to the question in greater detail with conditional characteristic perspective. The conditional characteristic function evolution, a new perspective in the Kalman filtering setting, is not available yet in literature. Secondly, the Kalman filtering can be sketched from the stochastic evolution of the conditional characteristic function for the linear Stochastic Differential Equation (SDE) with linear observation equation as well.

*Notation*: For the action of the conditional expectation operator, we adopt three notations, i.e. $\widehat{x}_t = E(x_t|Y_t) = \langle x_t \rangle$. A choice of the notation hinges on the structure of conditional expectation evolution equations.

**The Kalman filter: a conditional characteristics function perspective**

In this paper, we wish to construct the Kalman filter from the conditional characteristic function perspective by exploiting the Itô differential rule for the Wiener process [9], which is not available yet in literature. This paper is inspired from the fact that the older, but the celebrated result with new perspectives, produces greater insights. In this paper, the Itô calculus [10] and conditional characteristic function [5] are the cornerstone formalisms of the construction of the coupled Kalman filtering equations. Thus, we take a little pause and explain very briefly the influence of the Itô calculus and conditional characteristic function [2] in stochastic processes. The Itô calculus, a 1945 discovery, sowed the seed of stochastic calculus, a new branch of mathematics. The beauty, power and universality of the Itô calculus can be found in a famous paper of Kunita [10]. Later that the Itô calculus was adjudged for the inaugural Gauss prize. Secondly, the conditional characteristic function is a conditional moment generating function that forms the cornerstone formalism for the Central Limit Theorem of stochastic processes.

Consensus on the Itô-Stratonovich dilemma recommends 'the Stratonovich SDE is expected for stochastic processes in typical, continuous and real physical systems'. On the other hand, the Itô SDE is the correct stochastic calculus for stochastic finance and theoretical biology. The non-linear filtering would be different for two perspectives, the Itô and the Stratonovich. Interestingly, the Kalman filtering is invariant under the Itô and Stratonovich perspectives, see [11]. In this paper, we construct three Theorems, the first is about the *evolution* of conditional characteristic function for the Kalman filtering for the continuous state-discrete measurement system. The second Theorem is about the *stochastic evolution* of the conditional characteristic function [5] for the Kalman filtering for the continuous state-continuous measurement systems. The third Theorem is about the stochastic evolution of the bilinear Kalman filter for the continuous state-continuous measurement system, where the Kalman filtering is a special case of the bilinear Kalman filtering. This paper chiefly hinges on the evolution of conditional characteristic function and the evolution of conditional characteristic function holds for the continuous state. Thus, the continuous state obeying the Itô SDE is the subject of investigations in lieu of the stochastic difference equation.

**Theorem 1.** Consider the filtering model for the continuous state-discrete measurement system is

$$dx_i(t) = \sum_\alpha A_{i\alpha}(t)x_\alpha(t)dt + \sum G_{i\phi}(t)dW_\phi(t)$$

$$y_i(t_k) = \sum C_{i\phi}(t_k)x_\phi(t_k) + V_i(t_k),$$

where the vector Brownian motion process $W = \{W_t \in R^m, 0 \leq t < \infty\}$. Consider $x_t$ and $A_t$ denote the state vector and system state matrix respectively. The term $G(t)dW_t$ denotes the vector stochastic correction term with $\text{var}(dW_t) = Idt$. The term $y_{t_k}$ indicates a noisy observation vector. The observation noise $V_{t_k}$ is $N(0, R_k)$, where the term $R_k$ is a real-valued observation noise variance matrix. Then, the evolution of conditional characteristic function of the above filtering model between the observation instants, $\tau$, where $t_{k-1} \leq \tau < t_k$ becomes

$$d\langle \exp(s^T x_\tau) \rangle = \left\langle \sum_p \sum_\alpha A_{p\alpha}(\tau) x_\alpha \frac{\partial \exp(s^T x_\tau)}{\partial x_p} + \frac{1}{2}\sum_{p,q}(GG^T)_{pq}(\tau)\frac{\partial^2 \exp(s^T x_\tau)}{\partial x_p \partial x_q} \right\rangle d\tau.$$

Alternatively,

$$d\langle \exp(s^T x_\tau) \rangle = (\sum_p \sum_\alpha A_{p\alpha}(\tau) \left\langle x_\alpha(\tau)\frac{\partial \exp(s^T x_\tau)}{\partial x_p}\right\rangle$$

$$+ \frac{1}{2}\sum_{p,q}(GG^T)_{pq}(\tau)\left\langle\frac{\partial^2 \exp(s^T x_\tau)}{\partial x_p \partial x_q}\right\rangle)d\tau, \qquad (1)$$

where the characteristic function $\langle \exp(s^T x_\tau) \rangle = E(\exp(s^T x_\tau)|Y_{t_{k-1}})$. The conditional expectation operator $E$ accounts for the accumulated observation $Y_{t_{k-1}} = \{y_{t_n}, 0 \leq n \leq k-1\}$, $t_{k-1} < \tau \leq t_k$ and $t_{k-1}, t_k$ are the observation instants.

The Kalman filtering *between the observations* becomes a special case of the conditional characteristic function evolution. On the other hand, the Kalman filtering at the observation is not a consequence of the evolution of conditional characteristic function.

**Proof.** The proof of the conditional characteristic function evolution unaccounting for noisy observables can be sketched by exploiting the property of the differential and the expectation operators, Itô calculus. For the action of the conditional expectation operator $E(.)$ on lengthier terms, we adopt another notation $\langle \ \rangle$ as well.

$$\langle d \exp(s^T x_t) \rangle = \langle d\varphi(x_t, s) \rangle = \langle \varphi(x_{t+dt}, s) - \varphi(x_t, s) \rangle = \psi(t+dt, s) - \psi(t, s).$$

In the sense of characteristic function, the variable $s$ becomes $j\omega$, the input argument of the Fourier Transform. Now,

$$\langle \exp(s^T x_t) \rangle = \psi(t, s) = \langle \varphi(s, x_t) \rangle,$$

$$d\langle \exp(s^T x_t) \rangle = \psi(t+dt, s) - \psi(t, s) = \langle d \exp(s^T x_t) \rangle. \tag{2}$$

The above relation reveals the fact that the differential and the conditional expectation operator can be interchanged. Since we wish to develop the conditional characteristic function evolution at the instant $\tau$, where observations are not available, we adopt the notation $\tau$ in place of the time variable $t$. Thus,

$$d \exp(s^T x_\tau) = (\sum_p \sum_\alpha A_{p\alpha}(\tau) x_\alpha(\tau) \frac{\partial \exp(s^T x_\tau)}{\partial x_p} + \frac{1}{2} \sum_{p,q} (GG^T)_{pq}(\tau) \frac{\partial^2 \exp(s^T x_\tau)}{\partial x_p \partial x_q}) d\tau$$

$$+ \sum_{p,\phi} G_{p\phi}(\tau) \frac{\partial \exp(s^T x_\tau)}{\partial x_p} dW_\phi. \tag{3}$$

The above equation is a consequence of the Itô differential rule (Karatzas and Shreve [9]) for the scalar function $\exp(s^T x_t)$. From equation (2) and equation (3), we have

$$d\langle \exp(s^T x_t) \rangle = \left\langle \sum_p \sum_\alpha A_{p\alpha}(\tau) x_\alpha \frac{\partial \exp(s^T x_\tau)}{\partial x_p} + \frac{1}{2} \sum_{p,q} (GG^T)_{pq}(\tau) \frac{\partial^2 \exp(s^T x_\tau)}{\partial x_p \partial x_q} \right\rangle d\tau$$

$$+ \left\langle \sum_{p,\phi} G_{p\phi}(t) \frac{\partial \exp(s^T x_t)}{\partial x_p} dW_\phi \right\rangle \tag{4}$$

Since the conditional expectation operator is a linear operator and

$$\left\langle \sum_{p,\phi} G_{p\phi}(t) \frac{\partial \exp(s^T x_t)}{\partial x_p} dW_\phi \right\rangle = 0.$$

Equation (4) boils down to

$$d\langle \exp(s^T x_\tau) \rangle = \left\langle \sum_p \sum_\alpha A_{p\alpha}(\tau) x_\alpha(\tau) \frac{\partial \exp(s^T x_\tau)}{\partial x_p} + \frac{1}{2} \sum_{p,q} (GG^T)_{pq}(\tau) \frac{\partial^2 \exp(s^T x_\tau)}{\partial x_p \partial x_q} \right\rangle d\tau$$

$$= (\sum_p \sum_\alpha A_{p\alpha}(\tau) \left\langle x_\alpha(\tau) \frac{\partial \exp(s^T x_\tau)}{\partial x_p} \right\rangle + \frac{1}{2} \sum_{p,q} (GG^T)_{pq}(\tau) \left\langle \frac{\partial^2 \exp(s^T x_\tau)}{\partial x_p \partial x_q} \right\rangle) d\tau.$$

Thus, we arrive at equation (1) of the Theorem 1 of the paper.

*QED*

Now, consider $\exp(s^T x_t) = x_i(t)$. Then, the first term of the right-hand side will contribute and the last term of the right-hand side will vanish. As a result of this, equation (1) becomes

$$d\langle x_i(\tau)\rangle = \sum_\alpha A_{i\alpha}(\tau)\langle x_\alpha(\tau)\rangle d\tau,$$

Alternatively

$$d\langle x(\tau)\rangle = A_\tau \langle x_\tau\rangle d\tau. \tag{5}$$

Consider $\exp(s^T x_t) = x_i.x_j$, we have

$$d\langle x_i x_j\rangle = (\sum_\alpha A_{i\alpha}(\tau)\langle x_\alpha(\tau)x_j(\tau)\rangle + \sum_\alpha A_{j\alpha}(\tau)\langle x_\alpha(\tau)\rangle\langle x_j(\tau)\rangle$$

$$+\frac{1}{2}\sum_{p,q}(GG^T)_{pq}(\tau)\langle \frac{\partial^2 x_i x_j}{\partial x_p \partial x_q}\rangle)d\tau.$$

$$= (\sum_\alpha A_{i\alpha}(\tau)\langle x_\alpha(\tau)x_j(\tau)\rangle + \sum_\alpha A_{j\alpha}(\tau)\langle x_\alpha(\tau)x_i(\tau)\rangle$$

$$+\frac{1}{2}\sum_{p,q}(GG^T)_{pq}(\tau)\langle \delta_{ip}\delta_{jq} + \delta_{iq}\delta_{jp}\rangle)d\tau$$

$$= (\sum_\alpha A_{i\alpha}(\tau)\langle x_\alpha(\tau)x_j(\tau)\rangle + \sum_\alpha A_{j\alpha}(\tau)\langle x_\alpha(\tau)x_i(\tau)\rangle + (GG^T)_{ij}(\tau))d\tau.$$

After combining the above equation with equation (5), we get

$$dP_{ij} = d\langle x_i x_j\rangle - \langle x_i\rangle d\langle x_j\rangle - \langle x_j\rangle d\langle x_i\rangle - d\langle x_i\rangle d\langle x_j\rangle$$

$$= (\sum_\alpha A_{i\alpha}(\tau)\langle x_\alpha(\tau)x_j(\tau)\rangle + \sum_\alpha A_{j\alpha}(\tau)\langle x_\alpha(\tau)x_i(\tau)\rangle + (GG^T)_{ij}(\tau)$$

$$-\sum_\alpha A_{i\alpha}(\tau)\langle x_j(\tau)\rangle\langle x_\alpha(\tau)\rangle - \sum_\alpha A_{j\alpha}(\tau)\langle x_\alpha(\tau)\rangle\langle x_i(\tau)\rangle)d\tau$$

$$= (\sum_\alpha A_{i\alpha}(\tau)P_{j\alpha}(\tau) + \sum_\alpha P_{i\alpha}(\tau)A_{j\alpha}(\tau) + (GG^T)_{ij}(\tau))d\tau.$$

Alternatively,

$$dP_\tau = (A_\tau P_\tau + P_\tau A_\tau^T + (GG^T)_{ij}(\tau))d\tau. \tag{6}$$

Equation (1) describes the evolution of characteristic function for the Kalman filtering *between the observations*. Equations (5)-(6) are the coupled Kalman filtering equations *between the observations*. The Kalman filtering *at the observation* can be sketched using the conditional probability density and conditional expectation *at the observations* instants $t_{k-1}$ and $t_k$. A proof of the Kalman filtering equations *at the observation* can be found in a famous book authored by Andrew H Jazwinski on *Stochastic Processes* and *Non-linear filtering*, see [12] as well. That is not a consequence of the evolution of conditional characteristic function. Here, we omit the greater detail. For the vector case, we have equations (7)-(8).

$$\hat{x}_{t_k}^{t_k} = \hat{x}_{t_k}^{t_{k-1}} + P_{t_k}^{t_{k-1}} C_{t_k}^T (C_{t_k} P_{t_k}^{t_{k-1}} C_{t_k}^T + R_k)^{-1}(y_{t_k} - C_{t_k}\hat{x}_{t_k}^{t_{k-1}}), \tag{7}$$

$$P_{t_k}^{t_k} = P_{t_k}^{t_{k-1}} - P_{t_k}^{t_{k-1}} C_{t_k}^T (C_{t_k} P_{t_k}^{t_{k-1}} C_{t_k}^T + R_k)^{-1} C_{t_k} P_{t_k}^{t_{k-1}}. \tag{8}$$

Note that $\hat{x}_{t_k}^{t_k} = E(x_{t_k}|Y_{t_{k-1}}), \hat{x}_{t_k}^{t_k} = E(x_{t_k}|Y_{t_k})$. The Kalman filtering for the continuous state-discrete measurement system is given by the following set of four equations by adopting the vec function setting as well, i.e.

$$dvec(\langle x(\tau)\rangle) = vec(A_\tau \langle x_\tau \rangle)d\tau = (\langle x_\tau \rangle^T \otimes I)vec(A_\tau)d\tau, \quad (9a)$$

$$dvec(P_\tau) = ((I \otimes A_\tau + A_\tau \otimes I)vec(P_\tau) + vec((GG^T)(\tau)))d\tau, \quad (9b)$$

$$vec(\hat{x}_{t_k}^{t_k}) = vec(\hat{x}_{t_k}^{t_{k-1}}) + ((y_{t_k} - C_{t_k}\hat{x}_{t_k}^{t_{k-2}})^T (C_{t_k} P_{t_k}^{t_{k-1}} C_{t_k}^T + R_k)^{-1} C_{t_k} \otimes I) vec(P_{t_k}^{t_{k-1}}), \quad (9c)$$

$$vec(P_{t_k}^{t_k}) = vec(P_{t_k}^{t_{k-1}}) - (P_{t_k}^{t_{k-1}} C_{t_k}^T (C_{t_k} P_{t_k}^{t_{k-1}} C_{t_k}^T + R_k)^{-1} C_{t_k} \otimes I)vec(P_{t_k}^{t_{k-1}}). \quad (9d)$$

Equations (9a)-(9d) describe the Kalman filtering for the continuous state-discrete measurement system using the *vec* function interpretation. This paper is intended to develop the Kalman filtering in the conditional characteristic function evolution setting. On the other hand, the conditional characteristic function evolution setting is not applicable to the discrete state-discrete measurement system. Thus, in this paper, the Kalman filtering for the discrete state-discrete measurement system is not the subject of investigations. A good source on the Kalman filtering for the discrete state-discrete measurement can be found in Moore and Anderson [3].

**Remark 1**. The stationary Kalman filtering for the continuous state-discrete measurement has algebraic variance equations and constant coefficients. Conditional variance equations for *the between the observations case* and at *the observation case* become

$$0 = (I \otimes A + A \otimes I)vec(P_{t_{k-1}}^{t_{k-1}}) + vec(GG^T),$$

$$vec(P_{t_k}^{t_k}) = vec(P_{t_{k-1}}^{t_{k-1}}) - (P_{t_{k-1}}^{t_{k-1}} C_{t_k}^T (C_{t_k} P_{t_{k-1}}^{t_{k-1}} C_{t_k}^T + R_k)^{-1} C_{t_k} \otimes I)vec(P_{t_{k-1}}^{t_{k-1}})$$

respectively.

**Theorem 2**. Consider the *filtering model* for the continuous state-continuous measurement

$$dx_i(t) = \sum_\alpha A_{i\alpha}(t)x_\alpha(t)dt + \sum_\phi G_{i\phi}(t)dW_\phi(t)$$

$$dz_i(t) = \sum_\alpha C_{i\alpha}(t)x_\alpha(t)dt + d\eta_i,$$

where the vector Brownian motion process $W = \{W_t \in R^m, 0 \le t < \infty\}$. Consider $x_t$ and $A_t$ denote the state vector and system state matrix respectively. The term $G(t)dW_t$ denotes the vector stochastic correction term with $\text{var}(dW_t) = I dt$. The term $z_t$ indicates the noisy observation vector that is available continuously. The observation noise statistics is $N(0, \varphi_\eta dt)$. Then, the stochastic evolution of conditional characteristic function of the above filtering model assumes the structure of an SDE, i.e.

$$d\langle \exp(s^T x_t)\rangle = \left\langle \sum_p \sum_\alpha A_{p\alpha}(t)x_\alpha \frac{\partial \exp(s^T x_t)}{\partial x_p} + \frac{1}{2}\sum_{p,q}(GG^T)_{pq}(t) \frac{\partial^2 \exp(s^T x_t)}{\partial x_p \partial x_q} \right\rangle dt$$

$$+ (\langle \exp(s^T x_t)x_t^T C_t^T \rangle - \langle \exp(s^T x_t)\rangle\langle x_t^T C_t^T \rangle)\varphi_\eta^{-1}(dz_t - C_t \langle x_t \rangle dt).$$

$$= (\sum_p \sum_\alpha A_{p\alpha}(t)\left\langle x_\alpha(\tau)\frac{\partial \exp(s^T x_t)}{\partial x_p}\right\rangle + \frac{1}{2}\sum_{p,q}(GG^T)_{pq}(t)\left\langle \frac{\partial^2 \exp(s^T x_t)}{\partial x_p \partial x_q}\right\rangle)dt$$

$$+ \sum_{\alpha,\beta}(\left\langle \exp(s^T x_t)\sum_\gamma C_{\alpha\gamma}(t)x_\gamma(t)\right\rangle - \left\langle \exp(s^T x_t)\right\rangle\left\langle \sum_\gamma C_{\alpha\gamma}(t)x_\gamma(t)\right\rangle)$$

$$\times (\varphi_\eta^{-1})_{\alpha\beta} (dz_\beta - \sum_\gamma C_{\beta\gamma}(t)\langle x_\gamma \rangle dt), \qquad (10)$$

where the conditional characteristic function $\langle \exp(s^T x_t) \rangle = E(\exp(s^T x_t)|Y_t)$. The conditional expectation operator $E$ accounts for the accumulated observation $Y_t = \{z_\tau, t_0 \leq \tau \leq t\}$. The Kalman filtering becomes a special case of stochastic evolution of the conditional characteristic function, equation (10).

**Proof.** Here, we explain succinctly the proof of Theorem 2 of the paper. For the filtering model of Theorem 2, the stochastic evolution of conditional characteristic function assumes the structure of a stochastic integro-differential equation. One can construct the stochastic evolution of conditional characteristic function using Itô differential rule and the multiplication theorem of densities. As a result of this, we get equation (10) of Theorem 2 of the paper.

Notably, we can construct the stochastic evolution of conditional characteristic function by an alternative backward method. First, the method involves the construction of the stochastic evolution of conditional characteristic function for the Itô SDE using Itô differential rule and the multiplication theorem of densities (Liptser and Shiryaev [5]). As a result of this, the stochastic evolution of conditional characteristic function for the non-linear Itô SDE becomes

$$d\langle \exp(s^T x_t) \rangle = (\langle \sum_i f_i(t, x_t) \frac{\partial \exp(s^T x_t)}{\partial x_i} \rangle + \frac{1}{2} \langle \sum_{i,j} (GG^T)_{ij}(t, x_t) \frac{\partial^2 \exp(s^T x_t)}{\partial x_i \partial x_j} \rangle)dt$$
$$+ (\langle \varphi h^T \rangle - \langle \varphi \rangle\langle h^T \rangle)\varphi_\eta^{-1}(dz_t - \langle h \rangle dt), \qquad (11)$$

and the stochastic evolution for the linear Itô SDE is

$$d\langle \exp(s^T x_t) \rangle = (\sum_p \sum_\alpha A_{p\alpha}(t)\langle x_\alpha(\tau) \frac{\partial \exp(s^T x_t)}{\partial x_p} \rangle + \frac{1}{2} \sum_{p,q} (GG^T)_{pq}(t) \langle \frac{\partial^2 \exp(s^T x_t)}{\partial x_p \partial x_q} \rangle)dt$$
$$+ \sum_{\alpha,\beta} (\langle \exp(s^T x_t) \sum_\gamma C_{\alpha\gamma}(t)x_\gamma(t) \rangle - \langle \exp(s^T x_t) \rangle\langle \sum_\gamma C_{\alpha\gamma}(t)x_\gamma(t) \rangle)$$
$$\times (\varphi_\eta^{-1})_{\alpha\beta}(dz_\beta - \sum_\gamma C_{\beta\gamma}(t)\langle x_\gamma \rangle dt).$$

Thus, we arrive at equation (10) of Theorem 2 of the paper.

*QED*

That is a special case of equation (11) of the paper. Equation (10) can be regarded as the stochastic evolution of conditional characteristic function for the Kalman filtering for the continuous state-continuous measurement system. The term

$$(\langle \exp(s^T x_t)x_t^T C_t^T \rangle - \langle \exp(s^T x_t) \rangle\langle x_t^T C_t^T \rangle)\varphi_\eta^{-1}$$

of equation (10) can be regarded as the gain of the Kalman filter from the conditional characteristic function perspective.

Consider $\exp(s^T x_t) = x_i(t)$. Then, the first term of the right-hand side will contribute and the second term of the right-hand side will vanish. The last term of the above will contribute. Thus, equation (10) reduces to

$$d\langle x_i \rangle = (\sum_p \sum_\alpha A_{p\alpha}(t)\langle x_\alpha(\tau)\frac{\partial x_i}{\partial x_p}\rangle + \frac{1}{2}\sum_{p,q}(GG^T)_{pq}(t)\langle \frac{\partial^2 x_i}{\partial x_p \partial x_q}\rangle)dt$$

$$+ \sum_{\alpha,\beta}(\langle x_i \sum_\gamma C_{\alpha\gamma}(t)x_\gamma(t)\rangle - \langle x_i\rangle\langle \sum_\gamma C_{\alpha\gamma}(t)x_\gamma(t)\rangle)$$

$$\times (\varphi_\eta^{-1})_{\alpha\beta}(dz_\beta - \sum_\gamma C_{\beta\gamma}(t)\langle x_\gamma\rangle dt)$$

$$= (\sum_\alpha A_{i\alpha}(t)\langle x_\alpha(\tau)\rangle + \sum_{\alpha,\beta}(\langle \sum_\gamma C_{\alpha\gamma}(t)x_i(t)x_\gamma(t)\rangle - \langle x_i\rangle\langle \sum_\gamma C_{\alpha\gamma}(t)x_\gamma(t)\rangle)$$

$$\times (\varphi_\eta^{-1})_{\alpha\beta}(dz_\beta - \sum_\gamma C_{\beta\gamma}(t)\langle x_\gamma\rangle dt)$$

$$= (\sum_\alpha A_{i\alpha}(t)\langle x_\alpha(t)\rangle + \sum_{\alpha,\beta}(\sum_\gamma C_{\alpha\gamma}(t)(\langle x_i(t)x_\gamma(t)\rangle - \langle x_i\rangle\langle x_\gamma\rangle)$$

$$\times (\varphi_\eta^{-1})_{\alpha\beta}(dz_\beta - \sum_\gamma C_{\beta\gamma}(t)\langle x_\gamma\rangle dt)$$

$$= (\sum_\alpha A_{i\alpha}(t)\langle x_\alpha(t)\rangle dt + \sum_{\alpha,\beta}(\sum_\gamma C_{\alpha\gamma}(t)P_{i\gamma})(\varphi_\eta^{-1})_{\alpha\beta}(dz_\beta - \sum_\gamma C_{\beta\gamma}(t)\langle x_\gamma\rangle dt). \quad (12)$$

Consider $\exp(s^T x_t) = x_i.x_j$, equation (10) boils down to

$$d\langle x_i x_j\rangle = (\sum_\alpha A_{i\alpha}(t)\langle x_\alpha(t)x_j(t)\rangle + \sum_\alpha A_{j\alpha}(t)\langle x_\alpha(t)x_i(t)\rangle$$

$$+ \frac{1}{2}\sum_{p,q}(GG^T)_{pq}(t)\langle \frac{\partial^2 x_i x_j}{\partial x_p \partial x_q}\rangle)dt$$

$$+ \sum_{\alpha,\beta}(\langle x_i x_j \sum_\gamma C_{\alpha\gamma}(t)x_\gamma(t)\rangle - \langle x_i x_j\rangle\langle \sum_\gamma C_{\alpha\gamma}(t)x_\gamma(t)\rangle)$$

$$\times (\varphi_\eta^{-1})_{\alpha\beta}(dz_\beta - \sum_\gamma C_{\beta\gamma}(t)\langle x_\gamma\rangle dt). \quad (13)$$

After embedding equation (A.2) in equation (13), the above boils down to

$$(\sum_\alpha A_{i\alpha}(t)\langle x_\alpha(t)x_j(t)\rangle + \sum_\alpha A_{j\alpha}(t)\langle x_\alpha(t)x_i(t)\rangle$$

$$+ (GG^T)_{ij}(t))dt + \sum_{\alpha,\beta}\sum_\gamma (P_{i\gamma}\langle x_j\rangle + P_{j\gamma}\langle x_i\rangle)C_{\alpha\gamma}(t)$$

$$\times (\varphi_\eta^{-1})_{\alpha\beta}(dz_\beta - \sum_\gamma C_{\beta\gamma}(t)\langle x_\gamma\rangle dt),$$

$$\langle x_i(t)\rangle d\langle x_j(t)\rangle = (\sum_\alpha A_{j\alpha}(t)\langle x_i(t)\rangle\langle x_\alpha(t)\rangle$$

$$+ \sum_{\alpha,\beta}(\sum_\gamma C_{\alpha\gamma}(t)P_{j\gamma}\langle x_i(t)\rangle)(\varphi_\eta^{-1})_{\alpha\beta}(dz_\beta - \sum_\gamma C_{\beta\gamma}(t)\langle x_\gamma\rangle dt), \quad (14)$$

$$\langle x_j(t)\rangle d\langle x_i(t)\rangle = (\sum_\alpha A_{i\alpha}(t)\langle x_j(t)\rangle\langle x_\alpha(t)\rangle$$

$$+ \sum_{\alpha,\beta}(\sum_\gamma C_{\alpha\gamma}(t)P_{i\gamma})\langle x_j(t)\rangle(\varphi_\eta^{-1})_{\alpha\beta}(dz_\beta - \sum_\gamma C_{\beta\gamma}(t)\langle x_\gamma\rangle dt), \quad (15)$$

$$\begin{aligned}
d\langle x_i(t)\rangle d\langle x_j(t)\rangle &= (\sum_\alpha A_{i\alpha}(t)\langle x_\alpha(t)\rangle dt \\
&\quad + \sum_{\alpha,\beta}(\sum_\gamma C_{\alpha\gamma}(t)P_{i\gamma})\,(\varphi_\eta^{-1})_{\alpha\beta}(dz_\beta - \sum_\gamma C_{\beta\gamma}(t)\langle x_\gamma\rangle dt)) \\
&\quad \times (\sum_\alpha A_{j\alpha}(t)\langle x_\alpha(t)\rangle dt \\
&\quad + \sum_{\alpha,\beta}(\sum_\gamma C_{\alpha\gamma}(t)P_{j\gamma})\,(\varphi_\eta^{-1})_{\alpha\beta}(dz_\beta - \sum_\gamma C_{\beta\gamma}(t)\langle x_\gamma\rangle dt)) \\
&= (\sum_{\alpha,\beta}(\sum_\gamma C_{\alpha\gamma}(t)P_{i\gamma})\,(\varphi_\eta^{-1})_{\alpha\beta}(dz_\beta - \sum_\gamma C_{\beta\gamma}(t)\langle x_\gamma\rangle dt)) \\
&\quad \times (\sum_{r,s}(\sum_\gamma C_{s\gamma}(t)P_{j\gamma})\,(\varphi_\eta^{-1})_{rs}(dz_r - \sum_\gamma C_{r\gamma}(t)\langle x_\gamma\rangle dt)) \\
&= \sum_{(\alpha,\beta,r,s)}(\sum_\gamma C_{\alpha\gamma}(t)P_{i\gamma})(\sum_\gamma C_{s\gamma}(t)P_{j\gamma})(\varphi_\eta^{-1})_{\alpha\beta}(\varphi_\eta^{-1})_{rs}(\varphi_\eta)_{\beta r} \\
&= \sum_{(\alpha,\beta)}(\sum_\gamma C_{\alpha\gamma}(t)P_{i\gamma})(\sum_\gamma C_{\beta\gamma}(t)P_{j\gamma})(\varphi_\eta^{-1})_{\alpha\beta}(\varphi_\eta^{-1})_{\alpha\beta}(\varphi_\eta)_{\alpha\beta} \\
&= \sum_{(\alpha,\beta)}(\sum_\gamma C_{\alpha\gamma}(t)P_{i\gamma})(\sum_\gamma C_{\beta\gamma}(t)P_{j\gamma})(\varphi_\eta^{-1})_{\alpha\beta}. \qquad (16)
\end{aligned}$$

After combining equations (13)-(16), we have

$$\begin{aligned}
dP_{ij} &= d\langle x_i x_j\rangle - \langle x_i\rangle d\langle x_j\rangle - \langle x_j\rangle d\langle x_i\rangle - d\langle x_i\rangle d\langle x_j\rangle \\
&= (\sum_\alpha A_{i\alpha}(t)\langle x_\alpha(t)x_j(t)\rangle + \sum_\alpha A_{j\alpha}(t)\langle x_\alpha(t)x_i(t)\rangle \\
&\quad + (GG^T)_{ij}(t))dt + \sum_{\alpha,\beta}\sum_\gamma (P_{i\gamma}\langle x_j\rangle + P_{j\gamma}\langle x_i\rangle)C_{\alpha\gamma}(t) \\
&\quad \times (\varphi_\eta^{-1})_{\alpha\beta}(dz_\beta - \sum_\gamma C_{\beta\gamma}(t)\langle x_\gamma\rangle dt) \\
&\quad - (\sum_\alpha A_{j\alpha}(t)\langle x_i(t)\rangle\langle x_\alpha(t)\rangle)dt \\
&\quad - \sum_{\alpha,\beta}(\sum_\gamma C_{\alpha\gamma}(t)P_{j\gamma}\langle x_i(t)\rangle)\,(\varphi_\eta^{-1})_{\alpha\beta}(dz_\beta - \sum_\gamma C_{\beta\gamma}(t)\langle x_\gamma\rangle dt) \\
&\quad - \sum_\alpha A_{i\alpha}(t)\langle x_j(t)\rangle\langle x_\alpha(t)\rangle dt \\
&\quad - \sum_{\alpha,\beta}(\sum_\gamma C_{\alpha\gamma}(t)P_{i\gamma})\langle x_j(t)\rangle\,(\varphi_\eta^{-1})_{\alpha\beta}(dz_\beta - \sum_\gamma C_{\beta\gamma}(t)\langle x_\gamma\rangle dt) \\
&\quad - \sum_{(\alpha,\beta)}(\sum_\gamma C_{\alpha\gamma}(t)P_{i\gamma})(\sum_\gamma C_{\beta\gamma}(t)P_{j\gamma})(\varphi_\eta^{-1})_{\alpha\beta} dt \\
&= (\sum_\alpha A_{i\alpha}(t)(\langle x_\alpha(t)x_j(t)\rangle - \langle x_j(t)\rangle\langle x_\alpha(t)\rangle) \\
&\quad + \sum_\alpha A_{j\alpha}(t)(\langle x_\alpha(t)x_i(t)\rangle - \langle x_i(t)\rangle\langle x_\alpha(t)\rangle) + (GG^T)(t))dt \\
&\quad + \sum_{\alpha,\beta}\sum_\gamma (P_{i\gamma}\langle x_j\rangle + P_{j\gamma}\langle x_i\rangle)C_{\alpha\gamma}(t)\,(\varphi_\eta^{-1})_{\alpha\beta}(dz_\beta - \sum_\gamma C_{\beta\gamma}(t)\langle x_\gamma\rangle dt) \\
&\quad - \sum_{\alpha,\beta}(\sum_\gamma C_{\alpha\gamma}(t)P_{j\gamma}\langle x_i(t)\rangle)\,(\varphi_\eta^{-1})_{\alpha\beta}(dz_\beta - \sum_\gamma C_{\beta\gamma}(t)\langle x_\gamma\rangle dt)
\end{aligned}$$

$$-\sum_{\alpha,\beta}(\sum_{\gamma}C_{\alpha\gamma}(t)P_{i\gamma})\langle x_j(t)\rangle)\,(\varphi_\eta^{-1})_{\alpha\beta}(dz_\beta - \sum_{\gamma}C_{\beta\gamma}(t)\langle x_\gamma\rangle dt)$$

$$-\sum_{\alpha,\beta}(\sum_{\gamma}C_{\alpha\gamma}(t)P_{i\gamma})\,(\sum_{\gamma}C_{\beta\gamma}(t)P_{j\gamma})\,(\varphi_\eta^{-1})_{\alpha\beta}\,dt$$

$$= (\sum_{\alpha}A_{i\alpha}(t)P_{j\alpha}(t) + \sum_{\alpha}P_{i\alpha}(t)A_{j\alpha}(t) + (GG^T)_{ij}(t)$$

$$-\sum_{\alpha,\beta}(\sum_{\gamma}C_{\alpha\gamma}(t)P_{i\gamma})\,(\sum_{\gamma}C_{\beta\gamma}(t)P_{j\gamma})\,(\varphi_\eta^{-1})_{\alpha\beta})dt.$$

Alternatively, the standard Kalman filtering equations, the above equation coupled with equation (12), are recast as

$$d\langle x_t\rangle = A_t\langle x_t\rangle dt + P_t C_t^T \varphi_\eta^{-1}(dz_t - C_t\langle x_t\rangle dt),$$

$$dP_t = (A_t P_t + P_t A_t^T + (GG^T)(t) - P_t C_t^T \varphi_\eta^{-1} C_t P_t)dt.$$

The Kalman filtering for the continuous state-continuous state measurement system is given by the following set of two equations by adopting the vec function setting, i.e.

$$dvec(\langle x(t)\rangle) = (\langle x_t\rangle^T \otimes I)vec(A_t)dt + ((dz_t^T - \langle x_t\rangle^T C_t^T dt)\varphi_\eta^{-1}C_t \otimes I)vec(P_t), \quad (17a)$$

$$dvec(P_t) = ((I \otimes A_t + A_t \otimes I)vec(P_t) + vec((GG^T)(t))) - (P_t C_t^T \varphi_\eta^{-1} C_t \otimes I)vec(P_t). \quad (17b)$$

Equations (17a)-(178b) describe the Kalman filtering for the continuous state-continuous measurement system using the vec function interpretation.

**Remark 2.** The stationary Kalman filtering for the continuous state-continuous measurement has algebraic variance equation and constant coefficients. Alternatively, the conditional variance equation is

$$0 = (I \otimes A + A \otimes I)vec(P_t) + vec(GG^T) - (PC^T \varphi_\eta^{-1} C \otimes I)vec(P_t).$$

The following table 1 sums up five systems-theoretic properties of the Kalman filtering. Some of them are scattered in the literature and some of them are relatively very less known.

**Table 1**

| Properties | Kalman filtering for the continuous state- discrete measurement system | Kalman filtering for the continuous state- continuous measurement system |
|---|---|---|
| Evolution of conditional characteristic function [5] | Ordinary differential equation | Stochastic integro-differential equation |
| Riccati equation | Does not hold | Matrix ordinary differential equation |
| Lyapunov equation | Matrix ordinary differential equation | Does not hold |
| The Itô Vs Stratonovich [11] | Invariant | Invariant |
| Vec function [13] | $dvec(P_\tau) = ((I \otimes A_\tau + A_\tau \otimes I)vec(P_\tau) + vec((GG^T)(\tau)))d\tau,$ | $dvec(P_t) = ((I \otimes A_t + A_t \otimes I) \times vec(P_t)$ |

$$+ vec((GG^T)(t))$$
$$-(P_t C_t^T \varphi_\eta^{-1} C_t \otimes I) vec(P_t)) dt.$$

**Theorem 3.** Suppose the filtering model of a vector time-varying *bilinear* Stratonovich stochastic differential system is the following:
$$dx_t = (A_0(t) + A_t x_t) dt + (G(t) + x_t B_t^T) \circ dW_t,$$
$$dz_t = C_t x_t dt + d\eta_t,$$
where the solution $x_t \in U$, the phase space $U \subset R^n$. Note that the vector Brownian motion process $W = \{(W_1(t), W_2(t), ...., W_d(t)) \in R^d, 0 \le t < \infty\}$. The former part of the filtering model denotes the vector Brownian motion-driven vector time-varying Stratonovich bilinear stochastic differential equation. The latter denotes the noisy observation equation.

In the component-wise description, the coupled bilinear filtering equations for the filtering model are the following:

$$d\hat{x}_i(t) = (A_0^i(t) + \sum_\alpha A_{i\alpha}(t) \hat{x}_\alpha(t) + \frac{1}{2} \sum_\phi (G_{i\phi}(t) B_\phi(t) + B_\phi^2(t) \hat{x}_i(t))) dt$$
$$+ \sum_{\alpha,\beta} (\sum_p P_{ip}(t) C_{\alpha p}(t))(\varphi_\eta^{-1})_{\alpha\beta} (dz_\beta - \sum C_{\beta\gamma}(t) \hat{x}_\beta(t) dt), \quad (18a)$$

$$dP_{ij} = (\sum_p P_{ip} A_{jp}(t) + \sum_p P_{jp} A_{ip}(t) + \sum_\phi G_{i\phi}(t) G_{j\phi}(t) + \hat{x}_i \sum_\phi G_{j\phi}(t) B_\phi(t)$$
$$+ \sum_\phi G_{i\phi}(t) B_\phi(t) \hat{x}_j(t) + \hat{x}_i \hat{x}_j \sum_\phi B_\phi^2(t) + P_{ij} \sum_\phi B_\phi^2(t)$$
$$- \sum_{\alpha,\beta} (\sum_p P_{ip}(t) C_{\alpha p}(t)(\varphi_\eta^{-1})_{\alpha\beta} \sum_p P_{jp}(t) C_{\beta p}(t)) dt. \quad (18b)$$

The Kalman filtering becomes a special case of the coupled bilinear filtering equations, i.e. (18a)-(18b),

$$d\hat{x}_i(t) = (A_0^i(t) + \sum_\alpha A_{i\alpha}(t) \hat{x}_\alpha(t)) dt + \sum_{\alpha,\beta} (\sum_p P_{ip}(t) C_{\alpha p}(t))(\varphi_\eta^{-1})_{\alpha\beta} (dz_\beta - \sum C_{\beta\gamma}(t) \hat{x}_\beta(t) dt),$$

$$dP_{ij} = (\sum_p P_{ip} A_{jp}(t) + \sum_p P_{jp} A_{ip}(t) + \sum_\phi G_{i\phi}(t) G_{j\phi}(t) - \sum_{\alpha,\beta} (\sum_p P_{ip}(t) C_{\alpha p}(t))$$
$$- \sum_{\alpha,\beta} (\sum_p P_{ip}(t) C_{\alpha p}(t)) (\varphi_\eta^{-1})_{\alpha\beta} (\sum_p P_{jp}(t) C_{\beta p}(t))) dt.$$

**Proof.** The proof of the coupled bilinear filtering equations involves the following steps: First, construct the stochastic evolution of conditional characteristic function for the bilinear Stratonovich stochastic filtering model of Theorem 3 of the paper. Then, consider $\exp(s^T x_t) = x_i(t)$ and $\exp(s^T x_t) = x_i(t).x_j(t)$ to compute the component-wise bilinear stochastic filtering equations, equations (18a)–(18b). Since a procedure to construct the proof of the bilinear stochastic filtering equations is similar to the proof of the Theorem 2, we omit the detail. The following is worth to mention: (i) after ignoring the term,

$$\frac{1}{2} \sum_\phi (G_{i\phi}(t) B_\phi(t) + B_\phi^2(t) \hat{x}_i(t)) dt,$$

from the conditional mean equation (18a), we get the conditional mean equation of the Kalman filter. Secondly, after ignoring five correction terms from the conditional variance equation of the bilinear filter, equation (18b),

$$\widehat{x}_i \sum_{\phi} G_{j\phi}(t)B_{\phi}(t), \quad \sum_{\phi} G_{i\phi}(t)B_{\phi}(t)\widehat{x}_j(t),$$

$$\sum_{\phi} G_{i\phi}(t)B_{\phi}(t)\widehat{x}_j(t), \quad \widehat{x}_i\widehat{x}_j \sum_{\phi} B_{\phi}^2(t), \quad P_{ij} \sum_{\phi} B_{\phi}^2(t),$$

we arrive at the conditional variance evolution equation of the celebrated Kalman filter.

**Remark 3.** Theorem 3 is about the filtering equations for the bilinear Stratonovich stochastic differential equation. An alternative appealing interpretation is the bilinear Itô stochastic differential equation. The structure of the bilinear filtering equations hinges on stochastic interpretations. On the other hand, the Kalman filtering equations do not hinge on Itô and Stratonovich stochastic interpretations. The Itô-Stratonovich dilemma does not influence the Kalman filtering since the Itô and Stratonovich stochastic integrals coincide for linear stochastic differential equations.

**Conclusion**

In this paper, we have developed three Kalman filtering Theorems and their formal proof for the continuous state-discrete measurement system and continuous state-continuous measurement system. The paper exploits conditional characteristic function evolution equation and the Kiyoshi Itô calculus to achieve that. Most notably, this paper chiefly reveals 'how Kalman filtering equations become a special case of the stochastic evolution of conditional characteristic function for linear stochastic differential systems and linear measurement systems.

For the first time, the celebrated Kalman filtering was introduced in 1960. A philosophical contribution of the paper is to reveal connections between three influential results of stochastic processes, i.e. Kalman filtering, conditional characteristic function and the Itô calculus. This paper unifies them by adopting formal and systematic frameworks that were not available previously in the literature.

This paper differentiates between the Kalman filtering for the continuous state-continuous measurement system as well as the continuous state-discrete measurement system in table (1). Table (1) displays five major *ingredients* of the Kalman filter, conditional characteristic function evolution, the Riccati equation, the Lyapunov matrix equation, the vec function, the Itô-Stratonovich dilemma.

**Appendix**

Calculation of the term $\left\langle \exp(s^T x_t) \sum_{\gamma} C_{\alpha\gamma}(t)x_{\gamma}(t) \right\rangle - \left\langle \exp(s^T x_t) \right\rangle \left\langle \sum_{\gamma} C_{\alpha\gamma}(t)x_{\gamma}(t) \right\rangle$.

Here, we prove that

$$\left\langle \exp(s^T x_t) \sum_{\gamma} C_{\alpha\gamma}(t)x_{\gamma}(t) \right\rangle - \left\langle \exp(s^T x_t) \right\rangle \left\langle \sum_{\gamma} C_{\alpha\gamma}(t)x_{\gamma}(t) \right\rangle = \sum_{\gamma} C_{\alpha\gamma}(t)(P_{i\gamma}\left\langle x_j \right\rangle + P_{j\gamma}\left\langle x_i \right\rangle),$$

where the terms of the above relation are associated with the filtering model of the Theorem 2 of the paper.

The term $\left\langle \exp(s^T x_t) \sum_\gamma C_{\alpha\gamma}(t) x_\gamma(t) \right\rangle - \left\langle \exp(s^T x_t) \right\rangle \left\langle \sum_\gamma C_{\alpha\gamma}(t) x_\gamma(t) \right\rangle$ arises in the stochastic evolution of conditional characteristic function for the linear stochastic differential system coupled with the continuous linear measurement system. We wish to compute the term, where $\exp(s^T x_t) = x_i x_j$. Then, the resulting expression is embedded in the conditional variance evolution of the standard Kalman filter. Thus, we have

$$\left\langle \exp(s^T x_t) \sum_\gamma C_{\alpha\gamma}(t) x_\gamma(t) \right\rangle - \left\langle \exp(s^T x_t) \right\rangle \left\langle \sum_\gamma C_{\alpha\gamma}(t) x_\gamma(t) \right\rangle$$

$$= \left\langle x_i x_j \sum_\gamma C_{\alpha\gamma}(t) x_\gamma(t) \right\rangle - \left\langle x_i x_j \right\rangle \left\langle \sum_\gamma C_{\alpha\gamma}(t) x_\gamma(t) \right\rangle.$$

Since the conditional expectation operator is linear, the right-hand side of the above becomes

$$\left\langle \exp(s^T x_t) \sum_\gamma C_{\alpha\gamma}(t) x_\gamma(t) \right\rangle - \left\langle \exp(s^T x_t) \right\rangle \left\langle \sum_\gamma C_{\alpha\gamma}(t) x_\gamma(t) \right\rangle$$

$$= \left\langle x_i x_j \sum_\gamma C_{\alpha\gamma}(t) x_\gamma(t) \right\rangle - \left\langle x_i x_j \right\rangle \left\langle \sum_\gamma C_{\alpha\gamma}(t) x_\gamma(t) \right\rangle.$$

$$= \sum_\gamma C_{\alpha\gamma}(t) \left\langle x_i x_j x_\gamma \right\rangle - \sum_\gamma C_{\alpha\gamma}(t) \left\langle x_i x_j \right\rangle \left\langle x_\gamma \right\rangle. \quad (A.1)$$

After considering the Gaussian assumption, even powers will contribute and odd power contributions will vanish. Thus, the right-hand side of equation (A.1) becomes

$$\sum_\gamma C_{\alpha\gamma}(t) \left\langle x_i x_j x_\gamma \right\rangle - \sum_\gamma C_{\alpha\gamma}(t) \left\langle x_i x_j \right\rangle \left\langle x_\gamma \right\rangle$$

$$= \sum_\gamma C_{\alpha\gamma}(t) \left( \frac{1}{2} \sum_{p,q} P_{pq} \frac{\partial^2 \left\langle x_i \right\rangle \left\langle x_j \right\rangle \left\langle x_\gamma \right\rangle}{\partial \left\langle x_p \right\rangle \partial \left\langle x_q \right\rangle} - \frac{1}{2} \sum_{p,q} P_{pq} \frac{\partial^2 \left\langle x_i \right\rangle \left\langle x_j \right\rangle}{\partial \left\langle x_p \right\rangle \partial \left\langle x_q \right\rangle} \left\langle x_\gamma \right\rangle \right)$$

$$= \frac{1}{2} \sum_\gamma C_{\alpha\gamma}(t) \left( \sum_{p,q} P_{pq} \frac{\partial^2 \left\langle x_i \right\rangle \left\langle x_j \right\rangle \left\langle x_\gamma \right\rangle}{\partial \left\langle x_p \right\rangle \partial \left\langle x_q \right\rangle} - \sum_{p,q} P_{pq} \frac{\partial^2 \left\langle x_i \right\rangle \left\langle x_j \right\rangle}{\partial \left\langle x_p \right\rangle \partial \left\langle x_q \right\rangle} \left\langle x_\gamma \right\rangle \right)$$

$$= \frac{1}{2} \sum_\gamma C_{\alpha\gamma}(t) \sum_{p,q} P_{pq} \left( \frac{\partial}{\partial \left\langle x_p \right\rangle} (\delta_{iq} \left\langle x_j \right\rangle \left\langle x_\gamma \right\rangle + \delta_{jq} \left\langle x_i \right\rangle \left\langle x_\gamma \right\rangle + \delta_{q\gamma} \left\langle x_i \right\rangle \left\langle x_j \right\rangle) \right.$$

$$\left. - \delta_{ip} \delta_{jq} \left\langle x_\gamma \right\rangle - \delta_{jp} \delta_{iq} \left\langle x_\gamma \right\rangle \right).$$

The right-hand side of the above expression is further simplified to

$$\frac{1}{2} \sum_\gamma C_{\alpha\gamma}(t) \sum_{p,q} P_{pq} (\delta_{iq} \delta_{jp} \left\langle x_\gamma \right\rangle + \delta_{iq} \delta_{\gamma p} \left\langle x_j \right\rangle + \delta_{jq} \delta_{ip} \left\langle x_\gamma \right\rangle + \delta_{jq} \delta_{\gamma p} \left\langle x_i \right\rangle + \delta_{q\gamma} \delta_{ip} \left\langle x_j \right\rangle + \delta_{q\gamma} \delta_{jp} \left\langle x_i \right\rangle$$

$$- \delta_{ip} \delta_{jq} \left\langle x_\gamma \right\rangle - \delta_{jp} \delta_{iq}) \left\langle x_\gamma \right\rangle)$$

$$= \frac{1}{2}\sum_{\gamma} C_{\alpha\gamma}(t) \sum_{p,q} P_{pq} (\delta_{iq}\delta_{jp}\langle x_\gamma \rangle + \delta_{iq}\delta_{\gamma p}\langle x_j \rangle + \delta_{jq}\delta_{ip}\langle x_\gamma \rangle + \delta_{jq}\delta_{\gamma p}\langle x_i \rangle + \delta_{q\gamma}\delta_{ip}\langle x_j \rangle + \delta_{q\gamma}\delta_{jp}\langle x_i \rangle -$$
$$-\delta_{ip}\delta_{jq}\langle x_\gamma \rangle - \delta_{jp}\delta_{iq})\langle x_\gamma \rangle)$$
$$= \sum_{\gamma} C_{\alpha\gamma}(t) \sum_{p,q} P_{pq} (\delta_{iq}\delta_{\gamma p}\langle x_j \rangle + \delta_{jq}\delta_{\gamma p}\langle x_i \rangle)$$
$$= \sum_{\gamma} C_{\alpha\gamma}(t)(P_{i\gamma}\langle x_j \rangle + P_{j\gamma}\langle x_i \rangle).$$

Thus,

$$\left\langle \exp(s^T x_t)\sum_{\gamma} C_{\alpha\gamma}(t)x_\gamma(t) \right\rangle - \left\langle \exp(s^T x_t) \right\rangle \left\langle \sum_{\gamma} C_{\alpha\gamma}(t)x_\gamma(t) \right\rangle$$
$$= \left\langle x_i x_j \sum_{\gamma} C_{\alpha\gamma}(t)x_\gamma(t) \right\rangle - \left\langle x_i x_j \right\rangle \left\langle \sum_{\gamma} C_{\alpha\gamma}(t)x_\gamma(t) \right\rangle$$
$$= \sum_{\gamma} C_{\alpha\gamma}(t)(P_{i\gamma}\langle x_j \rangle + P_{j\gamma}\langle x_i \rangle). \quad (A.2)$$

*QED*

**REFERENCES**

1. Sharma, S. N., and Gawalwad, B. G., 2016, "Wiener meets Kolmogorov," *Norbert Wiener in the 21st Century* (Thinking Machines in the Physical World), Jul 13-Jul 15, University of Melbourne, Australia.
2. Kalman, R. E., 1960, "A New Approach to Linear Filtering and Prediction Problems," The ASME Transactions, Journal of Basic Engineering, 82, pp. 34-45.
3. Anderson, B. D. O., and Moore, J. B., 1991, "Kalman Filtering: Whence, What and Whither?," Mathematical System Theory, Ac Antoulas, ed., Springer, Berlin Heidelberg, pp. 41-54.
4. Gattami, A., 2014, "Kalman meets Shannon," *19th IFAC World Congress,* IFAC Proceedings, 47(3), pp. 2376–2381.
5. Liptser, R. S., and Shiryaev, A. N., 1977, *Statistics of Random Processes 1*, Springer, Berlin.
6. Fujisaki, M., Kallianpur, G., and Kunita, H., 1972, "Stochastic Differential Equations for the Nonlinear Filtering Problem," Osaka Journal of Mathematics, **9**(1), pp. 19-40.
7. Holbrook, J., Bhatia, R., 2006, "An Old Question Asked in a New Context Presents Strange Aspects, " Mathematical Intelligencer, 28 (1), pp. 32-39.
8. Anderson, B. D. O., and Moore, J. B., 1991, "Kalman Filtering: Whence, What and Whither? Mathematical System Theory ", Ac Antoulas, ed., Springer, Berlin Heidelberg, pp. 41-54.
9. Karatzas, I., and Shreve, S. E., 1991, *Brownian Motion and Stochastic Calculus*, Springer-Verlag, New York, Berlin.
10. Kunita, H., 2010, "Itô's Stochastic Calculus: Its Surprising Power for Applications, Stochastic Processes and their Applications," 120, pp.622-652.

11. Mannella, R., and McClintock, P. V. E., 2012, "Itô versus Stratonovich: 30 Years Later," Fluctuation and Noise Letters, **11**(1), pp.1240010-1240019.
12. Patel, H. G., and Sharma**,** S. N., 2014, "Third-order Continuous-discrete Filtering Equations for a Non-linear Dynamical System: The ASME Transactions, Journal of Computational and Non-linear Dynamics, **9**(3), pp. 034502-9.
13. Brewer, J.W., 1978, "Kronecker Products and Matrix Calculus in System Theory," IEEE Transactions on Circuits and Systems, **CAS-25 (**9), pp. 772-781.